\documentclass[twoside,12pt,a4paper]{amsart}
\usepackage{mathrsfs}
\usepackage{graphicx}
\usepackage{subfigure}
\usepackage{latexsym}
\usepackage{amsmath}
\usepackage{amsfonts}
\usepackage{amsthm}
\usepackage{amssymb}
\usepackage{color}
\usepackage{hyperref}
\hypersetup{colorlinks=true, citecolor=blue}

\newcommand{\etal}{\textit{et al.}}

\newtheorem{theo}{Theorem}[section]
\newtheorem{rem}{Remark}[section]

\newtheorem{lem}{Lemma}[section]
\newtheorem{prop}{Proposition}[section]
\newtheorem{defi}{Definition}[section]
\numberwithin{equation}{section}

\numberwithin{equation}{section}

\numberwithin{equation}{section}

\begin{document}
\title{\bf Structural stability of interior subsonic steady-states to  hydrodynamic model for semiconductors with sonic boundary}

\author{Yue-Hong Feng$^{1,4\text{\textdied}}$, Haifeng Hu$^{2,4\text{\textdied}}$ and Ming Mei$^{3,4\text{\textborn}}$}

\date{}

\maketitle\markboth{Y.~H.~Feng, H.~Hu and  M.~Mei } {Structural
stability of interior subsonic solutions with sonic boundary
conditions}

\vspace{-3mm}
\begin{center}
{\scriptsize
$^1$College of Mathematics, Faculty of Science, Beijing University of Technology, Beijing 100022, China \\[1mm]
$^2$School of Science, Changchun University, Changchun 130022, China\\[1mm]
$^3$Department of Mathematics, Champlain College Saint-Lambert, Quebec, J4P 3P2, Canada \\[1mm]
$^4$Department of Mathematics and Statistics, McGill University, Montreal, Quebec, H3A 2K6, Canada \\[1mm]
\textborn Corresponding author. E-mail: ming.mei@mcgill.ca\\[1mm]
Contributing authors. E-mails:  fyh@bjut.edu.cn; \hspace{1mm}  huhf@ccu.edu.cn \hspace{1mm}\\[1mm]
\textdied These authors contributed equally to this work.}
\end{center}

\
\begin{center}
\begin{minipage}{14cm}
{\bf Abstract.}{\small
For the stationary hydrodynamic model for semiconductors with sonic boundary, represented by Euler-Poisson equations, it possesses the various physical solutions including interior subsonic solutions/interior supersonic solutions/shock transonic solutions/$C^1$-smooth transonic solutions. However, the structural stability for these physical solutions is challenging and has remained open as we know.  In this paper, we  investigate the structural stability of interior subsonic solutions when the doping profiles are restricted in the  subsonic region. The main result is proved by using the local (weighted) singularity analysis and the monotonicity argument. Both the result itself and techniques developed here will give us some truly enlightening insights into our follow-up study on the structural stability of the remaining types of solutions.}
\end{minipage}
\end{center}

\vspace{7mm}

\vspace{5mm}

\vspace{3mm}
\baselineskip=15pt
\section{Introduction}
The hydrodynamic model was first derived by Bl{\o}tekj{\ae}r \cite{Bl70} for electrons in a semiconductor. After appropriate simplification the one-dimensional time-dependent system in the  isentropic case reads:
\begin{equation}\label{1.1}
\left\{
\begin{split}
&n_t+(nu)_x=0,\\
&(nu)_t+\left(nu^2+p(n)\right)_x=nE-\frac{nu}{\tau },\\
&E_x=n-b(x),
\end{split}
\right.
\end{equation}
where $n(x,t)$, $u(x,t)$ and $E(x,t)$ denote the electron density, velocity, and electric field respectively. The given function $p=p(n)$ is the pressure-density relation on which a commonly used hypothesis is
\begin{equation*}
p(n)=Tn^\gamma,
\end{equation*}
where $T>0$ is Boltzmann's constant and $\gamma\ge1$ is the adiabatic exponent. The constant parameter $\tau>0$ is the momentum relaxation time. The given background density $b(x)>0$ is called the doping profile. The hydrodynamic model \eqref{1.1} is also called Euler-Poisson equations with semiconductor effect. For more details we refer  to treatises \cite{MRS90,ZH16} and references therein.

In this paper, the main focus is on the isothermal steady-state flows satisfying equations
\begin{equation}\label{hu17}
\left\{
\begin{split}
&J\equiv\text{constant},\\
&\left(\frac{J^2}{n}+p(n)\right)_x=nE-\frac{J}{\tau},\\
&E_x=n-b(x),
\end{split}
\right.
\end{equation}
where $J=nu$ stands for the current density, and $p(n)=Tn$ corresponds to the isothermal ansatz. By the terminology from gas dynamics, we call $c:=\sqrt{P'(n)}=\sqrt T>0$ the speed of sound. The flow is referred to as subsonic, sonic or supersonic provided the velocity satisfies
\begin{equation}\label{1.5}
u<c,\quad u=c\quad\text{or}\quad u>c,\quad\text{respectively}.
\end{equation}
For convenience of notation, we introduce
\begin{equation}\label{1.2*}
\alpha=\frac{1}{\tau},~~\mbox{the reciprocal of the momentum relaxation time}.
\end{equation}
Without loss of generality, we set
\begin{equation}\label{z}
T=1 \mbox{ and } J=1,
\end{equation}
 thus the system \eqref{hu17} is equivalently reduced to the system
\begin{equation}\label{1.6}
\left\{
\begin{split}
&\left(1-\frac{1}{n^2}\right)n_x=nE-\alpha,\\
&E_x=n-b(x).
\end{split}
\right.
\end{equation}
From \eqref{1.5} and \eqref{z}, it is easy to see that the ﬂow is subsonic if $n>1$, sonic if $n=1$, or supersonic if $0<n<1$. By virtue of \eqref{1.2*}, we call the system \eqref{1.6} the Euler-Poisson equations with the semiconductor effect if $\alpha>0$, and without the semiconductor effect if $\alpha=0$, respectively. Throughout this paper, we are interested in the system \eqref{1.6} in the open interval $(0,1)$, which is subjected to the sonic boundary condition
\begin{equation}\label{1.7}
n(0)=n(1)=1.
\end{equation}
We also assume that the doping profile $b(x)$ is of class $C[0,1]$, satisfying the subsonic condition $b(x)>1$ on $[0,1]$. For simplicity of notation, its infimum and supremum over $[0,1]$ is denoted by
\[
\underline{b}:=\inf_{x\in[0,1]}b(x) \mbox{ and } \bar{b}:=\sup_{x\in[0,1]}b(x),
\]
 respectively.

Over the past three decades, major advances in the mathematical theory of steady-state Euler-Poisson equations with/without the semiconductor effect have been made by many authors. In what follows, we just list several results which are closely linked to the present paper.

For the purely subsonic steady-state flows,  in 1990, Degond
\etal\cite{DM90} first proved the existence of the subsonic solution
to the one-dimensional steady-state Euler-Poisson with the
semiconductor effect when its boundary states belongs to the
subsonic region. Subsequently, Degond \etal\cite{DM93} further
showed the existence and local uniqueness of irrotational subsonic
flows to the three-dimensional steady-state semiconductor
hydrodynamic model under a smallness assumptions on the data. Along
this line of research, the steady-state subsonic flows with and
without the semiconductor effect were investigated in various
physical boundary conditions and different dimensions
\cite{FI97,GS05,NS07,BDX14}. As for the purely supersonic
steady-state flows, Peng \etal\cite{PV06} established the existence
and uniqueness of the supersonic solutions with the semiconductor
effect, which correspond to a large current density.

Note that the system \eqref{hu17} or \eqref{1.6} will be degenerate
at the sonic state, thus the study on the transonic solutions and
various steady states satisfying the sonic boundary condition
becomes very difficult. Ascher \etal\cite{AMPS91} first examined the
existence of the transonic solution to the one-dimensional
isentropic Euler-Poisson equations without and with the
semiconductor effect when the doping profile is a supersonic
constant, and then Rosini \cite{Ro05} extended this work to the
non-isentropic case by the analysis of phase plane. When the doping
profile is non-constant, Gamba \cite{Ga92,GM96} investigated the
one-dimensional and two-dimensional transonic solutions with shocks,
respectively. However, these transonic solutions yield boundary
layers because they are constructed as the limits of vanishing
viscosity. Luo \etal\cite{LX12,LRXX11} further considered the
one-dimensional Euler-Poisson equations without the semiconductor
effect, under the restriction that boundary data are far from the
sonic state and the doping profile is either a subsonic constant or
a supersonic constant, a comprehensive analysis on the structure and
classification of steady states was carried out in \cite{LX12} by
using the analysis of phase plane. Meanwhile, both structural and
dynamical stability of steady transonic shock solutions was obtained
in \cite{LRXX11}.

What if the sonic state appears in the solutions? As we have seen,
all the existing works introduced above cannot answer this question.
Even the works regarding transonic shocks cannot radically answer it
either because the two different phase states are connected by the
jump of shocks satisfying the  Rankine-Hugoniot condition and
entropy condition, avoiding the degeneracy caused by the sonic
state. Recently, Li \etal\cite{LMZZ17,LMZZ18} systematically
explored the critical case, that is, the one-dimensional
semiconductor Euler-Poisson equations with the sonic boundary
condition. The existence, nonexistence and classification of all
types of physical steady states to this critical boundary-value
problem was obtained for the subsonic doping profile in
\cite{LMZZ17} and supersonic doping profile in \cite{LMZZ18}. More
precisely, in \cite{LMZZ17}, the authors proved that the critical
boundary-value problem admits a unique subsonic solution, at least
one supersonic solution, infinitely many transonic shocks if
$\alpha\ll1$, and infinitely many transonic $C^1$-smooth solutions
if $\alpha\gg1$; in \cite{LMZZ18}, the authors showed the
nonexistence of all types of physical steady states to the critical
boundary-value problem assuming that the doping profile is small
enough and $\alpha\gg1$, and they also discussed the existence of
supersonic and transonic shock solutions under the hypothesis that
the doping profile is close to the sonic state and $\alpha\ll1$.
Inspired by the groundbreaking works \cite{LMZZ17,LMZZ18}, there is
a series of interesting generalizations into the transonic doping
profile case in \cite{CMZZ20}, the case of transonic
$C^\infty$-smooth steady states in \cite{WMZZ21}, the
multi-dimensional cases in \cite{CMZZ21,CMZZ22}, and even the
bipolar case \cite{MMZ20}.

Compared with the existence theory in the critical case,
there are rather few works on the study of stability (both structural
 and dynamical stability) in the critical case due to the interior
  or boundary degeneracy. To the best of our knowledge,
  Feng \etal\cite{FMZ22} first demonstrated the structural
  stability of the $C^1$-smooth transonic steady states with
  respect to the small perturbation of both the supersonic
  constant doping profiles and non-degenerate boundary data.
  In \cite{FMZ22}, for supersonic doping profile, the authors
   also discussed the structural stability and linear dynamic
   instability of the transonic steady states with the non-degenerate
   boundary data under some appropriate hypotheses. Thus far, the most
   difficult part of the critical case, namely, the problem about the
    structural and dynamical stability of all types of physical solutions
     with the sonic boundary condition for the non-constant subsonic
     doping profile (for the existence theory, see \cite{LMZZ17}),
     is still open. To thoroughly solve this problem is full of challenges,
     owing to the boundary degeneracy. However, in the present paper,
     we intend to shed new light on this problem.

The purpose of this paper is to show that interior subsonic
solutions to the system \eqref{1.6} with the sonic boundary
condition \eqref{1.7} are structurally stable if we propose
a monotonicity restriction on subsonic doping profiles (see Theorem \ref{T2}).

This paper is organized as follows. Some necessary preliminaries
and the main result are stated in Section \ref{s2}. The proof
of the main result, Theorem \ref{T2}, is given in Section \ref{s3}.

\section{Preliminaries and the main result}\label{s2}
In this section we shall present the main result.
Before proceeding, we first give the important preliminaries
from the foregoing research \cite{LMZZ17}. First of all,
we recall the definition of the interior subsonic solution.
\begin{defi}\label{def}
We say a pair of functions $(n,E)(x)$ is an interior subsonic
solution of the boundary value problem \eqref{1.6}\&\eqref{1.7}
provided \textup{(i)} $(n-1)^2\in H_0^1(0,1)$, \textup{(ii)} $n(x)>1$,
for all $x\in(0,1)$, \textup{(iii)} $n(0)=n(1)=1$, \textup{(iv)}
the following equality holds for all test functions $\varphi\in H_0^1(0,1)$,
\begin{equation}\label{1.9}
\int_0^1\left(\frac{1}{n}-\frac{1}{n^3}\right)n_x\varphi_xdx+\alpha\int_0^1\frac{\varphi _x}{n}dx+\int_0^1\left(n-b\right)\varphi dx=0,
\end{equation}
and \textup{(v)} $E(x)$ is given by
\begin{equation}\label{E}
E(x)=\alpha +\int_0^x\left(n(y)-b(y)\right)dy.
\end{equation}
\end{defi}
In addition, we continue to recall the existence and uniqueness
of interior subsonic solutions, which is excerpted from the first
part of Theorem 1.3 in \cite{LMZZ17}.

\begin{prop}[Existence theory in \cite{LMZZ17}]\label{T1}
Suppose that the doping profile $b\in L^\infty(0,1)$ is subsonic
such that $\underline{b}>1$. Then for any $\alpha\in[0,\infty)$
the boundary value problem \eqref{1.6}\&\eqref{1.7} admits a unique
interior subsonic solution $(n,E)\in C^{\frac{1}{2}}[0,1]\times H^1(0, 1)$
satisfying the boundedness
\begin{equation}\label{ulb}
1+m\sin(\pi x)\le n(x)\le\bar{b}, \quad x\in [0,1],
\end{equation}
and the boundary behavior at endpoints
\begin{equation}\label{E(1)}
E(0)=\alpha,\quad E(1)<\alpha,
\end{equation}
\begin{equation}\label{L14}
\left\{
\begin{split}
& C_1(1-x)^\frac{1}{2}\le n(x)-1\le C_2(1-x)^\frac{1}{2},\\
& -C_3(1-x)^{-\frac{1}{2}}\le n_x(x)\le-C_4(1-x)^{-\frac{1}{2}},
\end{split}
\right. \quad\text{for}\ x\ \text{near}\ 1,
\end{equation}
where $m=m(\alpha, \underline{b})>0$, $C_2>C_1>0$ and $C_3>C_4>0$ are
certain uniform estimate constants.
\end{prop}

\begin{rem}
Note that the degeneracy of the boundary value problem
\eqref{1.6}\&\eqref{1.7} occurs merely on the boundary.
Thus, if we assume that the doping profile has relatively
higher-order regularity, say $b\in C[0,1]$, then by virtue
of the standard theory for elliptic interior regularity and
Sobolev's embedding theorem, the corresponding interior
subsonic solution $(n,E)$ is actually of class
$\left(C^1(0,1)\cap C^{\frac{1}{2}}[0,1]\right) \times C^1[0,1]$.
This fact will be tacitly exploited hereafter.
\end{rem}

We are now in a position to formulate the main result in the present paper.
\begin{theo}[Structural stability of interior subsonic solutions]\label{T2}
Assume that $b_1, b_2\in C[0,1]$ are subsonic such that $b_1(x)\ge b_2(x)>1$
for all $x\in[0,1]$. For $i=1,2$, let $(n_i,E_i)(x)$ denote the
interior subsonic solution relative to the doping profile $b_i$,
respectively. Then the two interior subsonic solutions are
structurally stable to one another in the sense that
\begin{equation}\label{1.10**}
\|n_1-n_2\|_{C[0,1]}+\|(1-x)^{\frac{1}{2}}(n_1-n_2)_x\|_{C[0,1]}+\|E_1-E_2\|_{C^1[0,1]}\le C\|b_1-b_2\|_{C[0,1]},
\end{equation}
where $C>0$ is a certain constant independent of $\|b_1-b_2\|_{C[0,1]}$.
\end{theo}

We conclude this section with a brief sketch of the strategy that
underlies the proof of our main result. Due to the boundary
degeneracy of interior subsonic solutions, the study of their
globally structural stability over the entire interval $[0,1]$
becomes sophisticated and challenging. Therefore, we shall have
to divide the whole interval $[0,1]$ into three domains as follows:
\begin{equation*}
[0,1]=[0,\delta)\cup[\delta,1-\delta]\cup(1-\delta,\delta],
\end{equation*}
where the intrinsic segmentation constant $\delta>0$ would be
appropriately determined (see Lemma \ref{L5}), and we will also
have to establish structural stability estimates separately on
their respective domains in the following order: (i) near the
left endpoint $x=0$; (ii) near the right endpoint $x=1$; (iii)
on the middle domain. This strategy is feasible because we have
discovered the following facts:
\begin{enumerate}
\item the local singularity analysis reveals that the plausible
singularity at the left endpoint $x=0$ is removable (see Lemma \ref{L1});
based on this, we are able to establish the local structural stability
 estimate on an inherent neighborhood $[0,\delta_0)$ by the
 monotonicity argument. The main point is that both the radius
 $\delta_0$ and the positive estimate constant are independent
 of $\|b_1-b_2\|_{C[0,1]}$ (see Lemma \ref{L3}). This is the
 reason why this type of neighborhood is referred to as ``to be intrinsic
 or inherent''. This sort of tacit convention will be used throughout
 the present paper.

\item the local weighted singularity analysis discloses that the
genuine singularity at the right endpoint $x=1$ can be well controlled
by the $(1-x)^{\frac{1}{2}}$-weight (see Lemma \ref{L4-1}); thus,
the monotonicity argument further ensures that the local weighted
structural stability holds on an intrinsic neighborhood
$(1-\delta_1,1]$ (see Lemma \ref{L4}).

\item the remaining part constitutes the middle domain, which is
regular as to the structural stability (see Lemma \ref{L5}).
\end{enumerate}

It is worth mentioning that the monotonicity argument has been playing
a crucial role in establishing structural stability estimates near
both endpoints. The principle behind the monotonicity argument
is given by Lemma \ref{Lcp}. From a technical point of view,
the monotonicity argument is useful, but only at the cost of
adding an extra restriction $b_1(x)\ge b_2(x)$ on $[0,1]$.
How to get rid of this restriction is a tough question, and
we will explore it in the future study.

\section{Proof of Theorem \ref{T2}}\label{s3}
This section is devoted to proving our main result. In order to make
the line of reasoning accessible to the reader,
the proof will be
divided into a sequence of lemmas.

We let $(n_i,E_i)(x)$ denote the interior subsonic solution
corresponding to the subsonic doping profile $b_i(x)>1$, satisfying
the sonic boundary value problem
\begin{equation}\label{1}
\left\{
\begin{split}
&\left(1-\frac{1}{n_i^2}\right)n_{ix}=n_iE_i-\alpha,\\
&E_{ix}=n_i-b_i(x), \quad x\in (0,1), \\
&n_i(0)=n_i(1)=1,
\end{split}
\right. \quad\text{for}\ i=1,2, \ \text{respectively}.
\end{equation}

First of all, we adapt the comparison principle in
\cite{LMZZ17}(Lemma 2.2, P4773) for use with two doping profiles
and their corresponding interior subsonic solutions, which is the
basis of the monotonicity argument in studying the structural
stability near endpoints.

\vspace{1.5mm}
\begin{lem}[Comparison principle]\label{Lcp}
Let the doping profiles $b_1, b_2\in C[0,1]$.
If $b_1(x)\ge b_2(x)>1$ on $[0,1]$. Then
\begin{equation}\label{1.10*}
n_1(x)\ge n_2(x),\ \text{on}\ [0,1].
\end{equation}
\end{lem}

\vspace{3mm}
\noindent\textbf{Proof.}
According to the relevant arguments from \cite{LMZZ17}(Equation (17), P4773),
for $i=1,2$, since $(n_i,E_i)$ is the interior subsonic solution,
thereby having the approximate solution sequence
$\{n_{ij}\}_{0<j<1}\subset C^1[0,1]$ satisfying the weak form
\begin{equation}\label{wf}
\int_0^1 A(n_{ij}, n_{ijx})\varphi_xdx+\int_0^1 (n_{ij}-b_i)\varphi dx=0,\quad\forall\varphi\in H_0^1(0,1),
\end{equation}
where
\begin{equation*}\label{1.12}
A(z,p):=\left(\frac{1}{z}-\frac{j^2}{z^3}\right)p+\alpha\frac{j}{z}.
\end{equation*}

Subtracting $\eqref{wf}|_{i=1}$ from $\eqref{wf}|_{i=2}$,
for all nonnegative test functions $\varphi\in H_0^1(0,1)$, we have
\begin{equation}\label{1.13}
\int_0^1 \left(A(n_{2j},n_{2jx})-A(n_{1j},n_{1jx})\right)\varphi_xdx+\int_0^1 \left(n_{2j}-n_{1j}\right)\varphi dx=\int_0^1 \left(b_2-b_1\right)\varphi dx\le0,
\end{equation}
where we have used the assumption that $b_1(x)\ge b_2(x)$ on $[0,1]$
in the last inequality. This is exactly the crucial Equation (19) in
\cite{LMZZ17}, the same result therefore applies to \eqref{1.13}
provided we simply imitate the remaining arguments in Lemma 2.2
of \cite{LMZZ17}. That is,
\begin{equation}\label{1.14}
n_{1j}(x)\ge n_{2j}(x),\ \text{on}\ [0,1], \ \text{for}\ 0<j<1.
\end{equation}

Now the monotonicity relation \eqref{Lcp} follows after a passage to
the limit as $j\to1^-$ on both sides of the inequality \eqref{1.14}.\qed
\vspace{3mm}

In addition, for $i=1,2$, we set about analyzing the boundary behavior of
the first-order derivative of $n_{i}(x)$ at the left endpoint $x=0$.
 It seems plausible that the singularity should have appeared there,
 as a matter of fact this ``fake'' singularity at $x=0$ is removable
 because of $E_i(0)=\alpha$.

\vspace{1.5mm}
\begin{lem}\label{L1}
Suppose that $b_i, i=1,2$ satisfy the same conditions in Lemma \ref{Lcp},
and $\alpha\ge2\sqrt{2}\max\{\sqrt{b_1(0)-1}, \sqrt{b_2(0)-1}\}$. Then
\begin{equation}\label{3}
\lim_{x\to0^+}n_{ix}(x)=\frac{1}{4}\left(\alpha-\sqrt{\alpha^2-8\left(b_i(0)-1\right)}\right)=:A_i>0,\quad i=1,2.
\end{equation}
\end{lem}

\vspace{3mm}
\noindent\textbf{Proof.}
In much the same way as in \cite{LMZZ17}(Theorem 5.6, P4802),
owing to $n_i(0)=1$ and $E_i(0)=\alpha$, it is easy to see
that $\lim_{x\to0^+}n_{ix}(x)$ exists by the monotone convergence
argument. Then from the first equation of \eqref{1}, we have
\begin{equation*}
n _{ix}=\frac{ E_i n _i^2 }{ n _i  + 1 } + \frac{\left(E_i-\alpha\right) n _i^2 }{ \left(  n _i  - 1  \right)\left( n _i  + 1  \right)}, \quad\text{in}\ (0,1).
\end{equation*}

Noting that $n_i(0)=1$ and $E_i(0)=\alpha$, it follows from the
L'Hospital Rule that
\begin{align*}
A_i &=\lim_{x\to0^+}n_{ix}(x)=\lim_{x\to0^+}\frac{E_in_i^2}{n_i+1}+\lim_{x\to0^+}\frac{\left(E_i-\alpha\right)n_i^2}{\left(n_i-1\right)\left(n_i+1\right)} \\
&=\frac{\alpha}{2}+\frac{1}{2}\lim_{x\to0^+}\frac{\left(E_i-\alpha\right)_x}{\left(n_i-1\right)_x}=\frac{\alpha}{2}+\frac{1}{2}\lim_{x\to0^+}\frac{E_{ix}}{n_{ix}}\\
&=\frac{\alpha}{2}+\frac{1}{2}\lim_{x\to0^+}\frac{n_i(x)-b_i(x)}{n_{ix}}\\
&=\frac{\alpha}{2}+\frac{1-b_i(0)}{2A_i},
\end{align*}
which in turn implies that
\begin{equation*}
A_i=\frac{1}{4}\left(\alpha-\sqrt{\alpha^2-8\left(b_i(0)-1\right)}\right)=\frac{2\left(b_i(0)-1\right)}{\alpha+\sqrt{\alpha^2-8\left(b_i(0)-1\right)}}=O(\frac{1}{\alpha}),
\end{equation*}
or
\begin{equation*}
A_i=\frac{1}{4}\left(\alpha+\sqrt{\alpha^2-8\left(b_i(0)-1\right)}\right)=O(\alpha).
\end{equation*}

According to the local singularity analysis in \cite{LMZZ17}(Lemma 5.3, P4796),
we know that in a small neighborhood of $x=0$, the drastic change of the
density component $n_i(x)$ of the interior subsonic solution is impossible
when $\alpha$ is suitably large. Therefore, we have to choose the former
root as the limit value of $\lim_{x\to0^+}n_{ix}(x)$, and the latter one is
the extraneous root.\qed
\vspace{3mm}

Based on Proposition \ref{T1} and Lemmas \ref{Lcp}$\sim$\ref{L1}, we are
now preparing to establish the local structural stability of interior
subsonic solutions to the boundary value problem \eqref{1.6}\&\eqref{1.7}
 on an intrinsic neighborhood of the left endpoint $x=0$.

\vspace{1.5mm}
\begin{lem}[Local structural stability estimate near $x=0$]\label{L3}
Under the same conditions in Lemma \ref{L1}. There exist two positive
constants $\delta_0\in(0,\frac{1}{2})$ and $C>0$ independent of
 $\|b_1-b_2\|_{C[0, 1]}$ such that
\begin{equation}\label{key-1}
\|n_1-n_2\|_{C^1[0,\delta_0)}+\|E_1-E_2\|_{C^1[0,\delta_0)}\le C\| b_1-b_2\|_{C[0,1]}.
\end{equation}
\end{lem}

\vspace{3mm}
\noindent\textbf{Proof.}
Firstly, in light of Lemma \ref{Lcp}, it is clear that the following
monotonicity relation holds,
\begin{equation}\label{mr}
\frac{n_1^3}{n_1+1}\ge\frac{n_2^3}{n_2+1},\quad\forall x\in[0,1].
\end{equation}

Next, for simplicity, we set $\tilde{E_i}:=E_i-\frac{\alpha}{n_i}$.
Multiplying Equation $\eqref{1}_1$ by $\frac{n_i^2}{n_i^2-1}$, we have
\begin{equation}\label{h1}
n_{ix}=\frac{\tilde{E}_i n_i^3}{n_i^2-1}, \quad  i=1,2.
\end{equation}
Taking the difference of Equations $\eqref{h1}|_{i=1}$ and $\eqref{h1}|_{i=2}$,
near $x=0$, we compute together with the monotonicity relation \eqref{mr} that
\begin{align}\label{33}
(n_1-n_2)_x=&\frac{\tilde{E}_1 n_1^3}{n_1^2-1}-\frac{\tilde{E}_2 n_2^3}{n_2^2-1}\\
=&\frac{n_1^3}{n_1+1}\frac{\tilde{E}_1}{n_1-1}-\frac{n_2^3}{n_2+1}\frac{\tilde{E}_1}{n_1-1}+\frac{n_2^3}{n_2+1}\frac{\tilde{E}_1}{n_1-1}-\frac{n_2^3}{n_2+1}\frac{\tilde{E}_2}{n_2 - 1}\notag\\
=&\frac{\tilde{E}_1}{n_1-1}\left(\frac{n_1^3}{n_1+1}-\frac{n_2^3}{n_2+1}\right)+\frac{n_2^3}{n_2+1}\left(\frac{\tilde{E}_1}{n_1-1}-\frac{\tilde{E}_2}{n_2-1}\right)\notag\\
\leq&M_0\alpha\left(\frac{n_1^3}{n_1+1}-\frac{n_2^3}{n_2+1}\right)+\frac{n_2^3}{n_2+1}M_0\|b_1-b_2\|_{C[0, 1]},\notag\\
\leq&C(n_1-n_2)+C\|b_1-b_2\|_{C[0, 1]},\quad x\in[0,\delta_0), \notag
\end{align}
where we have used the fact that there exist two positive constants
$\delta_0\in(0,\frac{1}{2})$ and $M_0>0$ independent of $\|b_1-b_2\|_{C[0, 1]}$ such that
\begin{equation}\label{0227}
\frac{\tilde{E}_1}{n_1- 1}(x)\leq M_0\alpha,\quad\text{and}\quad \left( \frac{\tilde{E}_1}{n_1- 1}-\frac{\tilde{E}_2}{n_2- 1}\right) (x)\leq M_0\|b_1-b_2\|_{C[0, 1]},\quad x\in[0,\delta_0).
\end{equation}

To prove that the crucial estimate \eqref{0227} on a certain intrinsic
neighborhood $[0,\delta_0)$ holds, we assume for the sake of contradiction
that for any $\delta\in(0,\frac{1}{2})$ and $M>0$, there exists
$x_\delta\in[0,\delta)$ such that
\begin{equation}\label{C0227}
\frac{\tilde{E}_1}{n_1-1}(x_\delta)>M\alpha, \quad\text{or}\quad \left( \frac{\tilde{E}_1}{n_1-1}-\frac{\tilde{E}_2}{n_2-1}\right) (x_\delta)> M\|b_1-b_2\|_{C[0, 1]}.
\end{equation}
Particularly, we take $\delta=\frac{1}{k}, k=3,4,5,\cdots$, for any $M>0$,
there is $x_k\in[0,\frac{1}{k})$ such that
\begin{equation*}
\frac{\tilde{E}_1}{n_1-1}(x_k)>M\alpha, \quad\text{or}\quad \left( \frac{\tilde{E}_1}{n_1-1}-\frac{\tilde{E}_2}{n_2-1}\right) (x_k)> M\|b_1-b_2\|_{C[0, 1]},
\end{equation*}
which implies that
\begin{equation}\label{0227-1}
\varliminf_{x_k\to0^+}\frac{\tilde{E}_1}{n_1-1}(x_k)\geq M\alpha,
\end{equation}
or
\begin{equation}\label{0227-2}
\varliminf_{x_k\to0^+}\left(\frac{\tilde{E}_1}{n_1-1}-\frac{\tilde{E}_2}{n_2-1}\right)(x_k)\geq M\|b_1-b_2\|_{C[0, 1]}.
\end{equation}
Combining the boundary behavior \eqref{E(1)}, the L'Hospital Rule,
Equation $\eqref{1}_2$ and Lemma \ref{L1}, we calculate
\begin{align}\label{0227-3}
\lim_{x\to0^+}\frac{\tilde{E}_i}{n_i-1}(x)=&\lim_{x\to0^+}\frac{E_i-\frac{\alpha}{n_i}}{n_i-1}(x)\\
=&\lim_{x\to0^+}\frac{E_i(x)-E_i(0)+\alpha-\frac{\alpha}{n_i(x)}}{n_i(x)-1} \notag\\
=&\lim_{x\to0^+}\frac{E_i(x)-E_i(0)}{n_i(x)-1}+\lim_{x\to0^+}\frac{\alpha}{n_i(x)}\notag\\
=&\lim_{x\to0^+}\frac{n_i-b_i}{n_{ix}}+\alpha \notag\\
=&\frac{1-b_i(0)}{A_i}+\alpha<\alpha, \quad i=1,2,\notag
\end{align}
and
\begin{align}\label{0227-4}
\lim_{x\to0^+}\left( \frac{\tilde{E}_1}{n_1-1}-\frac{\tilde{E}_2}{n_2-1}\right)(x)=&\frac{b_2(0)-1}{A_2}-\frac{b_1(0)-1}{A_1}\\
=&\frac{1}{2}\left(\sqrt{\alpha^2-8(b_2(0)-1)}-\sqrt{\alpha ^2-8(b_1(0)-1)}\right) \notag\\
=&\frac{2}{\sqrt{\alpha^2-8(\eta-1)}}\left(b_1(0)-b_2(0)\right) \notag\\
\leq&\tilde{C}_0\|b_1-b_2\|_{C[0, 1]},\notag
\end{align}
where $\eta\in(b_2(0), b_1(0))$. Furthermore,  we note that the constant $M$
in \eqref{0227-1} and \eqref{0227-2} can be chosen arbitrarily. Consequently,
if we take $M=2$ in \eqref{0227-1}, together with \eqref{0227-3}, we obtain
the contradiction that $2\alpha<\alpha$; if we take $M=2\tilde{C}_0$ in
\eqref{0227-2}, combined with \eqref{0227-4}, we have the contradiction $2\leq1$.

Based on the local estimate \eqref{33}, we continue establishing the structural
stability locally on the intrinsic neighborhood $[0,\delta_0)$. To this end,
we multiply through the inequality \eqref{33} by $n_1-n_2$ and calculate
\begin{equation}\label{26}
\frac{d}{dx}(  n_1 -n_2)^2(x)\le C( n_1-n_2)^2(x)+C\|b_1-b_2\|_{C[0,1]}^2,\quad x\in[0, \delta_0),
\end{equation}
where we have used Lemma \ref{Lcp} and Cauchy's inequality. By Gronwall's
inequality and the sonic boundary condition $n_1(0)=n_2(0)=1$, we get
\begin{equation}\label{27}
(n_1-n_2)^2(x)\le C\|b_1-b_2\|_{C[0, 1]}^2,\quad x\in[0,\delta_0),
\end{equation}
which in turn implies that
\begin{equation}\label{28}
|n_1-n_2|(x)+|\left(n_{1}-n_{2}\right)_x|(x)\le C\|b_1-b_2\|_{C[0, 1]},\quad x\in[0,\delta_0).
\end{equation}
with the aid of the foregoing local estimate \eqref{33} again.

Finally, from Equation \eqref{E} in Definition \ref{def}, we have
\begin{equation} \label{29}
E_i(x)=\alpha+\int_0^x\left(n_i(y)-b_i(y)\right)dy,\quad i=1,2.
\end{equation}
Taking the difference of $\eqref{29}|_{i=1}$ and $\eqref{29}|_{i=2}$, we
compute that
\begin{align}\label{30}
|E_1-E_2|(x)\le&\int_0^x|n_1-n_2|(y)dy+\int_0^x|b_1-b_2|(y)|dy\\
\le&C\|b_1-b_2\|_{C[0,1]},\quad x\in[0,\delta_0),\notag
\end{align}
and
\begin{align}\label{31}
|(E_1-E_2)_x|(x)=&|n_1-n_2-(b_1-b_2)|(x)\\
\le&C\|b_1-b_2\|_{C[0,1]}, \quad x\in[0,\delta_0).\notag
\end{align}
Hence, the local structural stability estimate \eqref{key-1} follows
immediately from Equations \eqref{28}, \eqref{30} and \eqref{31}.\qed
\vspace{3mm}

We now turn to analyzing the refined boundary behavior of the first-order
derivative of $n_i(x)$ at the right endpoint $x=1$. From the boundary estimate
 displayed in the second line of $\eqref{L14}$, we know that
 $\lim_{x\to1^-}n_{ix}(x)=-\infty$. This means the ``genuine''
 singularity will occur at the right endpoint $x=1$. Inspired
 by $\eqref{L14}$, we are able to implement the local ``weighted''
 singularity analysis. The result is summarized as follows.

\vspace{1.5mm}
\begin{lem}\label{L4-1}
Assume that $b_i, i=1,2$ satisfy the same conditions in Lemma \ref{Lcp}. Then
\begin{equation}\label{(1)}
\lim_{x\to1^-}(1-x)^{\frac{1}{2}}n_{ix}(x)=-\frac{1}{2}\sqrt{\int_0^1(b_i-n_i)dx}=:B_i<0,\quad i=1,2.
\end{equation}
\end{lem}

\vspace{3mm}
\noindent\textbf{Proof.}
For $i=1,2$, from the boundary estimate \eqref{L14}, we know that the
coefficient $1-\frac{1}{n_i^2}$ in the degenerate principal part of
Equation $\eqref{1}_1$ is comparable to $(1-x)^{\frac{1}{2}}$ near
the right endpoint $x=1$. Thus the regularity theory of
boundary-degenerate elliptic equations in one dimension
(e.g. \cite{Sc81}) ensures that $(1-x)^\frac{1}{2}n_{ix}(x)$ is
continuous up to the right endpoint $x=1$.

We now proceed to calculate the exact limit value of
$\lim_{x\to1-}(1-x)^{\frac{1}{2}}n_{ix}(x)$. For convenience, we set
\begin{equation*}
B_i:=\lim_{x\to1-}(1-x)^{\frac{1}{2}}n_{ix}(x).
\end{equation*}
Thereupon, multiplying through Equation $\eqref{1}_1$ by $(1-x)^{\frac{1}{2}}\frac{n_i^2}{n_i^2-1}$, we have
\begin{equation*}
(1-x)^\frac{1}{2}n_{ix}=\frac{n_i^3}{n_i+1}\left(E_i-\frac{\alpha}{n_i}\right)\frac{(1-x)^\frac{1}{2}}{n_i-1}.
\end{equation*}
By virtue of the sonic boundary condition $n_i(1)=1$, the known
boundary behavior \eqref{E(1)} and the L'Hospital Rule, we compute
\begin{align}\label{(3)}
B_i=&\lim_{x\to1^-}(1-x)^\frac{1}{2}n_{ix}\\
=&\lim_{x\to1^-}\frac{n_i^3}{n_i+1}\lim_{x\to1^-}(E_i-\frac{\alpha}{n_i})\lim_{x\to1^-}\frac{(1-x)^\frac{1}{2}}{n_i-1}\notag\\
=&\frac{1}{2}(E_i(1)-\alpha)\lim_{x\to1^-}\frac{-\frac{1}{2}(1-x)^{-\frac{1}{2}}}{n_{ix}}\notag\\
=&\frac{1}{4}(\alpha-E_i(1))\lim_{x\to1^-}\frac{1}{n_{ix}(1-x)^\frac{1}{2}}\notag\\
=&\frac{1}{4}(E_i(0)-E_i(1))\frac{1}{B_i},\notag
\end{align}
which implies from Equation $\eqref{1}_2$ that
\begin{equation*}
B_i=-\frac{1}{2}\sqrt{E_i(0)-E_i(1)}=-\frac{1}{2}\sqrt{\int_0^1(b_i-n_i)dx}<0,
\end{equation*}
where the boundary estimate $\eqref{L14}_2$ has been employed to uniquely
determine the value of $B_i$, which is strictly negative.\qed
\vspace{3mm}

Proposition \ref{T1} alongside Lemmas \ref{Lcp} and \ref{L4-1} now enable us
 to demonstrate the local weighted structural stability of interior subsonic
 solutions to the boundary value problem \eqref{1.6}\&\eqref{1.7} on an
 intrinsic neighborhood of the right endpoint $x=1$.

\vspace{1.5mm}
\begin{lem}[Local weighted structural stability estimate near $x=1$]\label{L4}
Under the same conditions in Lemma \ref{L1}. There exist two positive
constants $\delta_1\in(0,\frac{1}{2})$ and $C>0$ independent
of $\| b_1-b_2\|_{C[0, 1]}$ such that
\begin{align}\label{key-2}
&\left\|(1-x)^{-\frac{1}{2}}(n_1-n_2)\right\|_{C(1-\delta_1,1]}+\left\|(1-x)^\frac{1}{2}(n_{1x}-n_{2x})\right\|_{C(1-\delta_1,1]}\\
&+\|E_1-E_2\|_{C^1(1-\delta_1,1]}\le C\|b_1-b_2\|_{C[0,1]}.\notag
\end{align}
\end{lem}

\vspace{3mm}
\noindent\textbf{Proof.}
For $i=1,2$, from $\eqref{L14}_1$, we have known that $\frac{n_i-1}{(1-x)^{1/2}}$
 possesses the uniform positive upper and lower bounds near $x=1$, and
 so does its reciprocal $\frac{(1-x)^{1/2}}{n_i-1}$. This property will be
 used repeatedly hereafter.

Owing to the fact that $n_{ix}(x)$ has the genuine singularity at $x=1$, we
are compelled to establish the structural stability estimate near $x=1$ only
in the weighted manner as follows.

Firstly, multiplying through Equation $\eqref{1}_1$ by $(1-x)^{\frac{1}{2}}\frac{n_i^2}{n_i^2-1}$ and taking the difference of resultant equations for $i=1,2$, we calculate that
\begin{align}\label{(4)}
&(1-x)^\frac{1}{2}(n_{1x}-n_{2x})\\
=&\frac{n_1^2}{n_1+1}(n_1E_1-\alpha)\frac{(1-x)^\frac{1}{2}}{n_1-1}-\frac{n_2^2}{n_2+1}(n_2E_2-\alpha)\frac{(1-x)^\frac{1}{2}}{n_2-1}\notag\\
=&\frac{n_1^2}{n_1+1}(n_1E_1-\alpha)\left(\frac{(1-x)^\frac{1}{2}}{n_1-1}-\frac{(1-x)^\frac{1}{2}}{n_2-1}\right)\notag\\
&+\left(\frac{n_1^2}{n_1+1}(n_1E_1-\alpha)-\frac{n_2^2}{n_2+1}(n_2E_2-\alpha)\right)\frac{(1-x)^\frac{1}{2}}{n_2-1}\notag\\
=&h(n_1,E_1)\frac{(1-x)^\frac{1}{2}}{n_1-1}\frac{(1-x)^\frac{1}{2}}{n_2-1}\frac{n_2-n_1}{(1-x)^\frac{1}{2}}+\left(h(n_1,E_1)-h(n_2,E_2)\right)\frac{(1-x)^\frac{1}{2}}{n_2-1}\notag\\
=:&I_1+I_2,\notag
\end{align}
where
\begin{equation*}
h(n_i,E_i):=\frac{n_i^2}{n_i+1}(n_iE_i-\alpha),\quad i=1,2.
\end{equation*}

In what follows, near $x=1$, we shall estimate $I_1$ and $I_2$,
respectively. But first, we claim that the following estimates
\begin{gather}
|E_i(x)|\le\alpha+2\bar{b}_i,\quad x\in[0,1],\quad i=1,2,\label{hu1}\\
|E_1(1)-E_2(1)|\le C\|b_1-b_2\|_{C[0,1]}\label{hu2}
\end{gather}
hold, where the estimate constant $C>0$ is independent
of $\|b_1-b_2\|_{C[0,1]}$, and the proof of which is deferred
to Lemma \ref{lem-hu} at the end of this paper.

As for $I_1$, it is clear from $\eqref{L14}_1$ and \eqref{hu1} that
\begin{equation}\label{(5)}
|I_1|=\left|h(n_1,E_1)\frac{(1-x)^\frac{1}{2}}{n_1-1}\frac{(1-x)^\frac{1}{2}}{n_2-1}\frac{n_2-n_1}{(1-x)^\frac{1}{2}}\right|\le C\frac{|n_1-n_2|}{(1-x)^\frac{1}{2}}.
\end{equation}

However, as far as $I_2$ is concerned, the situation becomes more
complicated because of the factor $h(n_1,E_1)-h(n_2,E_2)$. Next, we are
taking it step by step. Precisely, a straightforward computation gives
\begin{align}\label{(6)}
h(n_1,E_1)-h(n_2,E_2)&=\frac{n_1^2}{n_1+1}(n_1E_1-\alpha)-\frac{n_2^2}{n_2+1}(n_2E_2-\alpha)\\
&=\alpha\left(\frac{n_2^2}{n_2+1}-\frac{n_1^2}{n_1+1}\right)+\left(\frac{n_1^3E_1}{n_1+1}-\frac{n_2^3E_2}{n_2+1}\right)\notag\\
&=:R_1+R_2.\notag
\end{align}
From \eqref{ulb} and Lemma \ref{Lcp}, we know that
\begin{equation}\label{hu0}
1\le1+m(\alpha,\underline{b}_2)\sin(\pi x)\le n_2(x)\le n_1(x)\le\overline{b}_1,\quad x\in[0,1].
\end{equation}
Consequently, it follows from the mean-value theorem of differentials that
\begin{equation}\label{(7)}
|R_1|=\left|\alpha\left(\frac{n_2^2}{n_2+1}-\frac{n_1^2}{n_1+1}\right)\right| \le C|n_1-n_2|\le C\frac{|n_1-n_2|(x)}{(1-x)^\frac{1}{2}}.
\end{equation}
We now turn to estimating $R_2$ near $x=1$. Combining \eqref{hu0},
\eqref{hu1}, \eqref{hu2}, the mean-value theorem of differentials,
and the mean-value theorem of integrals, we have
\begin{align}\label{hu3}
|R_2|&=\left|\frac{n_1^3E_1}{n_1+1}-\frac{n_2^3E_2}{n_2+1}\right|\\
&=\left|E_1\left(\frac{n_1^3}{n_1+1}-\frac{n_2^3}{n_2+1}\right)+\frac{n_2^3}{n_2+1}(E_1-E_2)\right|\notag\\
&\le C|n_1-n_2|(x)+C\left|\left[\big(E_1(1)-E_2(1)\big)-\left(\int_x^1(n_1-n_2)-(b_1-b_2)dy\right)\right]\right|\notag\\
&\le C|n_1-n_2|(x)+C|E_1(1)-E_2(1)|+C\int_x^1|n_1-n_2|dy+C\|b_1-b_2\|_{C[0,1]}\notag\\
&\le C\|b_1-b_2\|_{C[0,1]}+C\left(\frac{|n_1-n_2|(x)}{(1-x)^{\frac{1}{2}}}+\frac{|n_1-n_2|(\xi)}{(1-\xi)^{\frac{1}{2}}}\right),\quad\exists\xi\in[x,1],\notag
\end{align}
where we have used the formula
\begin{equation}
E_i(x)=E_i(1)-\int_x^1(n_i-b_i)(y)dy,\quad i=1,2.
\end{equation}
Substituting \eqref{(7)} and \eqref{hu3} into \eqref{(6)}, we have
\begin{equation*}
|h(n_1,E_1)-h(n_2,E_2)|\le C\|b_1-b_2\|_{C[0,1]}+C\left(\frac{|n_1-n_2|(x)}{(1-x)^{\frac{1}{2}}}+\frac{|n_1-n_2|(\xi)}{(1-\xi)^{\frac{1}{2}}}\right),\quad\exists\xi\in[x,1],
\end{equation*}
which further implies that
\begin{align}\label{hu4}
|I_2|&=\left|\left(h(n_1,E_1)-h(n_2,E_2)\right)\frac{(1-x)^\frac{1}{2}}{n_2-1}\right|\\
&\le C\|b_1-b_2\|_{C[0,1]}+C\left(\frac{|n_1-n_2|(x)}{(1-x)^{\frac{1}{2}}}+\frac{|n_1-n_2|(\xi)}{(1-\xi)^{\frac{1}{2}}}\right),\quad\exists\xi\in[x,1].\notag
\end{align}

Inserting \eqref{(5)} and \eqref{hu4} into \eqref{(4)}, near $x=1$, we
obtain \begin{align}\label{(14)}
&(1-x)^\frac{1}{2}|(n_1-n_2)_x|(x)\\
\le&C\left(\frac{|n_1-n_2|(x)}{(1-x)^{\frac{1}{2}}}+\frac{|n_1-n_2|(\xi)}{(1-\xi)^{\frac{1}{2}}}\right)+C\|b_1-b_2\|_{C[0,1]},\quad\exists\xi\in[x,1].\notag
\end{align}
It is worth mentioning that the generic constant $C>0$ in \eqref{(14)} is
independent of $\|b_1-b_2\|_{C[0,1]}$. Moreover,
the term $\frac{|n_1-n_2|(x)}{(1-x)^{1/2}}+\frac{|n_1-n_2|(\xi)}{(1-\xi)^{1/2}}$
on the right-hand side of \eqref{(14)} can be bounded by an appropriate
constant multiple of $\|b_1-b_2\|_{C[0,1]}$ in an intrinsic neighborhood
of the right endpoint $x=1$. Precisely, we claim that there exist two
positive constants $0<\delta_1<\frac{1}{2}$ and $M_1>0$ independent
of $\|b_1-b_2\|_{C[0,1]}$ such that
\begin{equation}\label{(14*)}
\frac{|n_1-n_2|(x)}{(1-x)^\frac{1}{2}}\le M_1\|b_1-b_2\|_{C[0,1]}, \quad x\in(1-\delta_1,1].
\end{equation}
Aiming for a contradiction, suppose that for any $\delta\in (0,\frac{1}{2})$
and $M>0$, there is $x_\delta\in(1-\delta,1]$ such that
\begin{equation}
\frac{|n_1-n_2|(x_\delta)}{(1-x_\delta)^\frac{1}{2}}>M\|b_1-b_2\|_{C[0,1]}.
\end{equation}
By the arbitrariness, we could take $\delta=\frac{1}{k}, k=3,4,5,\cdots$, for
arbitrary $M>0$, there exists $x_k\in(1-\frac{1}{k},1]$ such that
\begin{equation}
\frac{|n_1-n_2|(x_k)}{(1-x_k)^\frac{1}{2}}>M\|b_1-b_2\|_{C[0,1]},
\end{equation}
which implies that
\begin{equation}\label{0227-5}
\varliminf_{x_k\to1^-}\frac{|n_1-n_2|(x_k)}{(1-x_k)^\frac{1}{2}}\ge M\|b_1-b_2\|_{C[0,1]}.
\end{equation}
Besides, combining Lemma \ref{Lcp}, the L'Hospital Rule and Lemma
\ref{L4-1}, we calculate that
\begin{align}\label{0227-6}
\lim_{x\to1^-}\frac{|n_1-n_2|(x)}{(1-x)^\frac{1}{2}}=&\lim_{x\to1^-}\frac{(n_1-n_2)(x)}{(1-x)^\frac{1}{2}}=\lim_{x\to1^-}\frac{(n_1-1)(x)}{(1-x)^\frac{1}{2}}-\lim_{x\to1^-}\frac{(n_2-1)(x)}{(1-x)^\frac{1}{2}}\\
=&-2\lim_{x\to1^-}(1-x)^{\frac{1}{2}}n_{1x}+2\lim_{x\to1^-}(1-x)^{\frac{1}{2}}n_{2x}\notag\\
=&\sqrt{\int_0^1(b_1-n_1)dx}-\sqrt{\int_0^1(b_2-n_2)dx}\notag\\
=&\frac{\int_0^1(b_1-b_2)dx-\int_0^1(n_1-n_2)dx}{\sqrt{\int_0^1(b_1-n_1)dx}+\sqrt{\int_0^1(b_2-n_2)dx}}\notag\\
\le&\frac{\int_0^1(b_1-b_2)dx}{\sqrt{\int_0^1(b_1-n_1)dx}+\sqrt{\int_0^1(b_2-n_2)dx}}\notag\\
\le&\tilde{C}_1\|b_1-b_2\|_{C[0,1]}.\notag
\end{align}
Moreover, we note that the constant $M>0$ in \eqref{0227-5} is arbitrary.
Therefore, together with \eqref{0227-6}, taking $M=2\tilde{C}_1$
in \eqref{0227-5} leads to the contradiction that $2\le1$.

Applying \eqref{(14*)} to \eqref{(14)}, we have
\begin{equation}\label{hu5}
(1-x)^\frac{1}{2}|(n_1-n_2)_x|(x)\le C\|b_1-b_2\|_{C[0,1]}, \quad x\in(1-\delta_1,1].
\end{equation}

Similarly to \eqref{hu3}, we are able to compute that
\begin{align}\label{(21)}
|E_1-E_2|(x)\le&C\|b_1-b_2\|_{C[0,1]}+\frac{|n_1-n_2|(\xi)}{(1-\xi)^{\frac{1}{2}}},\quad\exists\xi\in[x,1]\\
\le&C\|b_1-b_2\|_{C[0,1]},\quad x\in(1-\delta_1,1],\notag
\end{align}
and
\begin{align}\label{(22)}
|(E_1-E_2)_x|(x)=&|n_1-n_2-(b_1-b_2)|(x)\le|n_1-n_2|(x)+|b_1-b_2|(x)\\
\le&\frac{|n_1-n_2|(x)}{(1-x)^{\frac{1}{2}}}+\|b_1-b_2\|_{C[0,1]}\notag\\
\le&C\|b_1-b_2\|_{C[0,1]},\quad x\in(1-\delta_1,1].\notag
\end{align}

Finally, putting results \eqref{(14*)}, \eqref{hu5}, \eqref{(21)}
and \eqref{(22)} together, we obtain the desired local weighted
estimate \eqref{key-2}.\qed
\vspace{3mm}

Up to now, we have obtained two intrinsic small domains $[0,\delta_0)$
and $(1-\delta_1,1]$ distributed around the two endpoints $x=0$ and $x=1$,
respectively. This fact enables us to establish the structural stability
estimate on a certain regular domain $[\delta,1-\delta]$, where
$0<\delta:=\min\{\delta_0,\delta_1\}<1/2$.

\vspace{1.5mm}
\begin{lem}\label{L5}
Under the same conditions in Lemma \ref{L1}. Let $\delta:=\min\{\delta_0,\delta_1\}$. Then there is a positive constant $C>0$ independent of $\|b_1-b_2\|_{C[0,1]}$ such that
\begin{equation}\label{key-3}
\|n_1-n_2\|_{C^1[\delta,1-\delta]}+\|E_1-E_2\|_{C^1[\delta,1-\delta]}\le C\|b_1-b_2\|_{C[0,1]}.
\end{equation}
\end{lem}

\vspace{3mm}
\noindent\textbf{Proof.}
We are now able to work with Equations \eqref{1} on a regular closed
interval $[\delta,1-\delta]$ away from singularities, where $\delta$
has been defined in the hypothesis of the present lemma.

Firstly, we rewrite the estimate \eqref{hu0} on the regular interval as follows.
\begin{equation}\label{hu6}
1<l:=1+m(\alpha,\underline{b}_2)\sin(\pi\delta)\le n_2(x)\le n_1(x)\le\overline{b}_1,\quad x\in[\delta,1-\delta].
\end{equation}

Secondly, subtracting $\eqref{1}|_{i=2}$ from $\eqref{1}|_{i=1}$, for $x\in[\delta,1-\delta]$, we thus get
\begin{align}\label{41}
(n_1-n_2)_x=&\frac{n_1^3E_1-\alpha n_1^2}{n_1^2-1}-\frac{n_2^3E_2-\alpha n_2^2}{n_2^2-1}\\
=&E_1\left(f(n_1)-f(n_2)\right)+f(n_2)(E_1-E_2)-\alpha\left(g(n_1)-g(n_2)\right)\notag\\
=&\left(E_1 f'(\bar{\eta})-\alpha g'(\tilde{\eta})\right)(n_1-n_2)+f(n_2)(E_1-E_2),\quad\exists\bar{\eta},\tilde{\eta}\in(n_2,n_1),\notag
\end{align}
and
\begin{equation}\label{42}
(E_1-E_2)_x =(n _1-n _2)-(b_1-b_2),
\end{equation}
where
\begin{equation*}
f(n):=\frac{n^3}{n^2-1},\quad g(n):=\frac{n^2}{n^2-1}, \quad\forall n\in[l,\overline{b}_1],
\end{equation*}
and we have used the mean-value theorem of differentials in the third
line of Equation \eqref{41}.

Thirdly, multiplying through \eqref{41} by $n_1-n_2$, and using \eqref{hu1},
\eqref{hu6} and Cauchy's inequality together, we have
\begin{equation}\label{hu7}
\left((n_1-n_2)^2\right)_x\le C(\alpha,l,\bar{b}_1)\left((n_1-n_2)^2+(E_1-E_2)^2\right), \quad x\in[\delta,1-\delta].
\end{equation}
Similarly, multiplying through \eqref{42} by $E_1-E_2$, and employing
Cauchy's inequality, we obtain
\begin{equation}\label{hu8}
\left((E_1-E_2)^2\right)_x\le(n_1-n_2)^2+2(E_1-E_2)^2+\|b_1-b_2\|_{C[0,1]}^2, \quad x\in[\delta,1-\delta].
\end{equation}
And then, summing estimates \eqref{hu7} and \eqref{hu8} gives
\begin{align}\label{hu9}
&\left((n_1-n_2)^2+(E_1-E_2)^2\right)_x(x)\\
\le&C\left((n_1-n_2)^2+(E_1-E_2)^2\right)(x)+\|b_1-b_2\|_{C[0,1]}^2, \quad x\in[\delta,1-\delta].\notag
\end{align}
Applying the Gronwall inequality to \eqref{hu9}, we have
\begin{align}\label{hu10}
&\left((n_1-n_2)^2+(E_1-E_2)^2\right)(x)\\
\le&e^{\int_\delta^xCdy}\left[\left((n_1-n_2)^2+(E_1-E_2)^2\right)(\delta)+\int_\delta^x\|b_1-b_2\|_{C[0,1]}^2dy\right]\notag\\
\le&C\left[\left((n_1-n_2)^2+(E_1-E_2)^2\right)(\delta)+\|b_1-b_2\|_{C[0,1]}^2\right], \quad x\in[\delta,1-\delta].\notag
\end{align}
Noting that $\delta\le\delta_0$ and the continuity of the error
function pair $(n_1-n_2, E_1-E_2)(x)$ at $x=\delta_0$, from
Lemma \ref{L3} we see
\begin{equation}\label{hu11}
\left((n_1-n_2)^2+(E_1-E_2)^2\right)(\delta)\le C\|b_1-b_2\|_{C[0,1]}^2,
\end{equation}
which along with \eqref{hu10} implies
\begin{equation}\label{47}
|n_1-n_2|(x)+|E_1-E_2|(x)\le C\|b_1-b_2\|_{C[0,1]}, \quad x\in[\delta,1-\delta].
\end{equation}

Finally, from Equations \eqref{41}\&\eqref{42}, and the estimate \eqref{47},
 we directly calculate
\begin{equation}\label{hu12}
|(n_1-n_2)_x|(x)+|(E_1-E_2)_x|(x)\le C\|b_1-b_2\|_{C[0,1]}, \quad x\in[\delta,1-\delta].
\end{equation}
Combining estimates \eqref{47} and \eqref{hu12} yields the desired structural
 stability estimate \eqref{key-3} on the regular domain $[\delta,1-\delta]$.\qed
\vspace{3mm}

Last but not least, let us prove the estimates \eqref{hu1} and \eqref{hu2} in
 the following lemma.

\vspace{1.5mm}
\begin{lem}\label{lem-hu}
Under the same conditions in Lemma \ref{L4}. Then there exists a
positive constant $C$ independent of $\|b_1-b_2\|_{C[0,1]}$ such that
 estimates \eqref{hu1} and \eqref{hu2} hold, that is,
\begin{equation*}
|E_i(x)|\le\alpha+2\bar{b}_i,\quad x\in[0,1],\quad i=1,2,
\end{equation*}
and
\begin{equation*}
|E_1(1)-E_2(1)|\le C\|b_1-b_2\|_{C[0,1]},
\end{equation*}
respectively.
\end{lem}

\vspace{3mm}
\noindent\textbf{Proof.}
From Equation \eqref{E} in Definition \ref{def}, we have
\begin{equation} \label{hu13}
E_i(x)=\alpha+\int_0^x(n_i-b_i)(y)dy,\quad\forall x\in[0,1], \quad i=1,2.
\end{equation}

First of all, in light of the lower and upper bounds \eqref{ulb} of $n_i(x)$,
a straightforward computation gives
\begin{align}
|E_i(x)|=\left|\alpha+\int_0^x(n_i-b_i)dy\right|\le\alpha+\int_0^1(n_i+b_i)dy\le\alpha+2\bar{b}_i,\quad\forall x\in[0,1], \quad i=1,2.
\end{align}

Next, taking the value $x=1$ in Equation \eqref{hu13}, we have
\begin{equation}\label{hu14}
E_i(1)=\alpha+\int_0^1(n_i-b_i)(y)dy, \quad i=1,2.
\end{equation}
Furthermore, taking the difference of Equations $\eqref{hu14}|_{i=1}$
and $\eqref{hu14}|_{i=2}$, we calculate
\begin{align}\label{hu15}
|E_1(1)-E_2(1)|\le&\left|\int_0^1\Big[(n_1-b_1)-(n_2-b_2)\Big]dy\right|\\
\le&\int_0^1|n_1-n_2|(y)dy+\|b_1-b_2\|_{C[0,1]}\notag\\
=&|n_1-n_2|(\xi)+\|b_1-b_2\|_{C[0,1]}, \quad\exists\xi\in[0,1],\notag
\end{align}
where we have used the mean-value theorem of integrals in the last line.

Finally, we claim that there is a positive constant $C$ independent of
$\|b_1-b_2\|_{C[0,1]}$ such that
\begin{equation}\label{hu16}
|n_1-n_2|(\xi)\le C\|b_1-b_2\|_{C[0,1]},
\end{equation}
wherever the point $\xi$ is located in the whole interval $[0,1]$.
In fact, take $\delta$ the same as in Lemma \ref{L5}, and if
$\xi\in[0,1-\delta]$, it is clear from estimates \eqref{key-1}
and \eqref{key-3} that \eqref{hu16} is true; if $\xi\in(1-\delta,1]$,
the intrinsic local estimate \eqref{(14*)} guarantees that \eqref{hu16}
is true as well. Consequently, substituting \eqref{hu16} into \eqref{hu15},
we obtain the desired estimate \eqref{hu2}.\qed
\vspace{3mm}

We end this section by a summary for the proof of Theorem \ref{T2}
because the previous lemmas have already made the proof evident.

\vspace{3mm}
\noindent\textbf{Proof of Theorem \ref{T2}.}
Obviously, putting all the estimates we have established in
Lemmas \ref{L3}, \ref{L5} and \ref{L4} together, we have the
globally structural stability estimate \eqref{1.10**}.\qed

\vspace{5mm}
\noindent{\sc Acknowledgments:} This work was commenced while the
first two authors were visiting McGill University from 2021 to 2022.
They would like to express their gratitude to McGill University for
its hospitality. The research of Y. H. Feng was supported by China
Scholarship Council for the senior visiting scholar program
(202006545001). The research of H. Hu was partially supported by
National Natural Science Foundation of China (Grant No.11801039),
China Scholarship Council (No.202007535001) and Scientific Research
Program of Changchun University (No.ZKP202013). The research of M.
Mei was partially supported by NSERC grant RGPIN 354724-2016.



\begin{thebibliography}{99}
\bibitem{AMPS91} U.M. Ascher, P.A. Markowich, P. Pietra and C. Schmeiser,
 A phase plane analysis of transonic solutions for the hydrodynamic
 semiconductor model, {\it Math. Models Methods Appl. Sci.},
 1(3) (1991), 347-376.

\bibitem{Bl70} K. Bl{\o}tekj{\ae}r, Transport equations for
electrons in two-valley semiconductors, {\it IEEE Trans.
Electron Devices}, 17 (1970), 38-47.

\bibitem{BDX14} M.  Bae,  B.  Duan  and  C.J.  Xie, Subsonic solutions
for steady Euler-Poisson system in two-dimensional nozzles,  {\it SIAM J.
Math. Anal.}, 46 (2014), 3455-3480.

\bibitem{CMZZ20} L. Chen, M. Mei, G. Zhang and K. Zhang, Steady
hydrodynamic model of semiconductors with sonic boundary and transonic
doping profile, {\it J. Differential Equations}, 269 (2020), 8173-8211.

\bibitem{CMZZ21} L. Chen, M. Mei, G. Zhang and K. Zhang, Radial solutions
 of the hydrodynamic model of semiconductors with sonic boundary,
  {\it J. Math. Anal. Appl.}, 501 (2021), 125187.

\bibitem{CMZZ22}L. Chen, M. Mei, G. Zhang and K. Zhang, Transonic
steady-states of Euler-Poisson equations for semiconductor models
with sonic boundary. {\it SIAM J. Math. Anal.}, 54 (2022), no. 1,
363-388.

\bibitem{DM90}  P. Degond and  P.A.  Markowich,  On  a one-dimensional
steady-state hydrodynamic model for semiconductors,
{\it Appl. Math. Lett.}, 3(3) (1990), 25-29.

\bibitem{DM93}  P. Degond and  P.A.  Markowich, A steady state
potential flow model for semiconductors, {\it  Ann. Mat. Pura Appl.},
165(4) (1993), 87-98.

\bibitem{FI97} W. Fang and K. Ito, Steady-state solutions of a
one-dimensional hydrodynamic model for semiconductors, {\it  J.
Differential Equations}, 133 (1997), 224-244.

\bibitem{FMZ22} Y.H.~Feng, M. Mei and G. Zhang, Nonlinear structural
stability and linear dynamic instability of transonic steady-states to
a hydrodynamic model for semiconductors, arXiv:2202.03475, submitted.

\bibitem{Ga92} I.M. Gamba, Stationary transonic solutions of a one-dimensional
hydrodynamic model for semiconductors.
{\it Commun. Partial Differ. Equ.}, 17(3-4) (1992), 553-577.

\bibitem{GM96} I.M. Gamba and C.S. Morawetz,  A viscous approximation
for a 2-D steady semiconductor or transonic gas dynamic flow: existence
 theorem for potential flow, {\it Commun. Pure Appl. Math.}, 49(10) (1996),
  999-1049.

\bibitem{GS05} Y. Guo and W. Strauss, Stability of semiconductor states
with insulating and contact boundary conditions,
{\it Arch. Rational Mech. Anal.}, 179 (2005),  1-30.

\bibitem{LMZZ17} J. Li, M. Mei, G. Zhang and K. Zhang, Steady
hydrodynamic model of semiconductors with sonic boundary: (I)
Subsonic doping profile, {\it SIAM J. Math. Anal.}, 49 (2017), 4767-4811.

\bibitem{LMZZ18} J. Li, M. Mei, G. Zhang and K. Zhang, Steady
hydrodynamic model of semiconductors with sonic boundary: (II)
Supersonic doping profile, {\it SIAM J. Math. Anal.}, 50 (2018),  718-734.

\bibitem{LRXX11} T. Luo, J. Rauch, C.J. Xie and Z.P. Xin, Stability of
transonic shock solutions for one-dimensional Euler-Poisson equations,
 {\it  Arch. Rational Mech. Anal.},  202 (2011), 787-827.

\bibitem{LX12} T. Luo and Z. Xin,  Transonic shock solutions for
a system of Euler-Poisson equations, {\it Commun. Math. Sci.}, 10
(2012), 419-462.

\bibitem{MRS90} P.A.~Markowich, C.A.~Ringhofer and C.~Schmeiser,
{\it Semiconductor equations}, Springer, 1990.

\bibitem{MMZ20} P. Mu, M. Mei and K. Zhang, Subsonic and supersonic
steady-states of bipolar hydrodynamic model of semiconductors with
sonic boundary. {\it Commun. Math. Sci.}, 18 (2020), no. 7,
2005-2038.

\bibitem{NS07} S. Nishibata and M. Suzuki, Asymptotic stability of a
stationary solution to a hydrodynamic model of semiconductors,
 {\it Osaka J. Math.}, 44 (2007), 639-665.

\bibitem{PV06} Y. Peng and I. Violet,  Example of supersonic solutions
to a steady state Euler-Poisson system, {\it Appl. Math. Lett.} 19(12)
(2006), 1335-1340.

\bibitem{Ro05} M. Rosini,  A phase analysis of transonic solutions for
the hydrodynamic semiconductor model. {\it Q. Appl. Math.}, 63(2) (2005),
251-268.

\bibitem{Sc81}
R. Schreiber, Regularity of singular two-point boundary value problems.
 {\it SIAM J. Math. Anal.}, 12 (1981), no. 1, 104-109.

\bibitem{WMZZ21} M. Wei, M. Mei, G. Zhang and K. Zhang, Smooth transonic
 steady-states of hydrodynamic model for semiconductors,
 {\it SIAM J. Math. Anal.}, 53(4), (2021), 4908-4932.

\bibitem{ZH16} K. Zhang and H. Hu, {\it Introduction to
Semiconductor Partial Differential Equations (Chinese Edition)},
Science Press, Beijing, 2016.
\end{thebibliography}
\end{document}